\documentclass[12pt]{article}

\usepackage[T2A]{fontenc}           % Русская статья
\usepackage[cp1251]{inputenc}       % Кодировка WINDOWS
\usepackage{amsfonts}
\usepackage[active]{srcltx}
\usepackage{amsmath}
\usepackage{amssymb}
\usepackage{cite}
\usepackage[dvips]{graphicx}
\usepackage{color}
\newenvironment{sciabstract}{%
\begin{quote} \footnotesize} %\bf
{\end{quote}}
\usepackage{graphicx}
\title{\bf To the numerical solution of the inverse multi-frequency scalar acoustics problem }
\author
{A.B. Bakushinsky$^1$, A.S. Leonov$^2$,\\\
\\
\normalsize{$^{1}$Federal Research Center “Computer Science and Control”}\\
\normalsize{of Russian Academy of Sciences,}
\normalsize{Institute for Systems Analysis}\\
\normalsize{Moscow, Russian Federation}\\
\normalsize{$^{2}$National Research Nuclear University MEPhI,}\\
\normalsize{Moscow, Russian Federation}\\
}
\date{}

\begin{document}

%%% remove comment delimiter ('%') and select language if required
%\selectlanguage{russian}

\maketitle

\begin{sciabstract}
A new algorithm is proposed for solving the three-dimensional scalar inverse problem of acoustic sounding in an inhomogeneous medium. The data for the algorithm are the complex amplitudes of the wave field measured outside the inhomogeneity region. For the data recording scheme “in a flat layer”, the inverse problem is reduced using the Fourier transform to solving a set of one-dimensional Fredholm integral equations of the first kind. Solving these equations by the use of regularizing methods, we calculate the complex amplitude of the wave field in the inhomogeneity region and then we find the desired field of sound velocities in this domain. The proposed algorithm allows us to solve the inverse problem on a personal computer of average performance (without parallelization) for sufficiently fine three-dimensional grids in a few minutes. We demonstrate the results of solving model inverse problems at one frequency and several frequencies simultaneously along with a study of the accuracy of the proposed algorithm. We also investigate the issues of numerical stability of the algorithm with respect to data perturbations.

\textbf{Pacs}: 02.30.Zz, 02.30.Jr, 02.60.Cb

\end{sciabstract}

\section{Introduction}
Let the scalar function $p(\textbf{x},t)$ determine the acoustic wave field depending on the coordinates $\textbf{x}\in {\bf {\mathbb R}}^{3} $ and time $t\ge 0$. The field is created in an infinite medium by sources localized in a known region, $ S $. The medium is characterized by the local phase velocity of sound $c(\textbf{x})$ and has a constant density. Moreover, it is known that $c(\textbf{x})=c_{0} =\mathrm{const}$ outside a given region $X$ satisfying the condition $X\cap S=\varnothing $. In the region $X$, the function $c(\textbf{x})$ can be variable, and this is interpreted as  acoustic inhomogeneities there. In this case, for a harmonic source of the form $f(\textbf{x},\omega )e^{i\omega t} $ with a known frequency $\omega $, the field $p(\textbf{x},t)$ can be found in the form $p(\textbf{x},t)=u(\textbf{x},\omega )e^{i\omega t} $ in the framework of linear acoustics, where the complex amplitude $u(\textbf{x},\omega )$ satisfies the equation
\begin{equation}\label{L1_1_}
\Delta u(\textbf{x},\omega )+k_{0}^2 u(\textbf{x},\omega )=f(\textbf{x},\omega )+\omega ^{2}
\xi (\textbf{x})u(\textbf{x},\omega ),\, \, \, \textbf{x}\in {\bf {\mathbb R}}^{3} ,
\end{equation}
and the Sommerfeld radiation condition (see, for example, \cite{1}). Here $k_{0} =\frac{\omega }{c_{0} } $ and $\xi (\textbf{x})=c_{0}^{-2} -c^{-2} (\textbf{x})$.
We are interested in the following inverse problem for the equation \eqref{L1_1_}. Knowing the complex amplitude of the field $u(\textbf{x},\omega )$ in the region $ Y, \,Y\cap X=\varnothing, \,Y\cap S=\varnothing  $, for a certain set of frequencies $\omega $, find the coefficient $\xi (\textbf{x})$, that is the function $c(\textbf{x})$, which determines the acoustic inhomogeneities in the domain $X$. Introducing the Green function for the Helmholtz equation \eqref{L1_1_} in ${\bf {\mathbb R}}^{3} $, $G(\rho ,\omega )=-\frac{\exp \left(i\omega \rho /c_{0} \right)}{4\pi \rho } $, we can reduce the inverse problem, under certain assumptions about the smoothness of the functions $u(\textbf{x},\omega ),f(\textbf{x},\omega ),\, c(\textbf{x})$, to non-linear system of integral equations with respect to unknowns $u(\textbf{x}',\omega ),\,\xi(\textbf{x}')\,(\textbf{x}'\in X)$ (see, e.g., \cite{1})
\begin{equation}\label{L1_2_}
\left\{ \begin{array}{l}
u({\bf{x}},\omega ) = {u_0}({\bf{x}},\omega ) + {\omega ^2}\int_X^{} {G(|{\bf{x}} - {\bf{x}}'|,\omega )\xi ({\bf{x}}')u({\bf{x}}',\omega )d{\bf{x}}'} ,\,\,\,{\bf{x}} \in X\\
\\
{\omega ^2}\int_X^{} {G(|{\bf{x}} - {\bf{x}}'|,\omega )\xi ({\bf{x}}')u({\bf{x}}',\omega )d{\bf{x}}' = W({\bf{x}},\omega ),\,\,{\bf{x}} \in Y.}
\end{array} \right.
\end{equation}
The functions included in \eqref{L1_2_},
\[u_{0} (\textbf{x},\omega )=\int _{X}^{}G(|\textbf{x}-\textbf{x}'|,\omega )f(\textbf{x}',\omega )d\textbf{x}' ,\, \, \, W(\textbf{x},\omega )=u(\textbf{x},\omega )-u_{0} (\textbf{x},\omega ),\, \, \textbf{x}\in Y,\]
are known (computable) functions, and the quantities
$u(\textbf{x}',\omega ),\, \xi (\textbf{x}'),\, \, \textbf{x}'\in X,$
must be determined.

The problem \eqref{L1_2_}, in particular, the existence and uniqueness of its solution, has been studied theoretically quite well (see, for example, \cite{21,22,2} and others).

In many papers, various numerical methods have been proposed for solving the problem \eqref{L1_2_} in two-dimensional and three-dimensional formulations. For example, in \cite{1}, the system \eqref{L1_2_} is reduced to a nonlinear operator equation, which is then solved by a special iterative method. However, in this approach, the authors do not take into account that this equation is an ill-posed problem and do not regularize their iterative method. Nevertheless, the method works for model “medium power” scatterers \cite[c.101-105]{1}.

In \cite{2}, the regularized Gauss – Newton method was used to solve the system of equations \eqref{L1_2_} in the three-dimensional axially symmetric case. In \cite{3}, a special gradient method and the Fletcher – Reeves method were used for a similar problem. Other gradient methods and a supercomputer were used in \cite{4,5}.

Note that there are alternative approaches to solving the inverse acoustic sounding problem that are not directly related to the system \eqref{L1_2_}. In particular, the original method of boundary control, which allows one to solve three-dimensional problems similar to \eqref{L1_2_}, was proposed and developed in \cite{6,7}. To solve the two-dimensional inverse scattering problem, the well-known method of R.G. Novikov \cite{9} was used in \cite{8}. In the subsequent work \cite{10}, a comparative analysis of a variant of this method and some other functional-analytical methods for solving two-dimensional inverse problems of acoustic scattering is performed. The methods of M.V. Klibanov, summarized in the monograph \cite{11}, and also the methods from the monograph \cite{Kab} turned out to be very promising when processing real experimental data. We also note recent works \cite{Kl1,Kl2}.

All the mentioned methods for solving the inverse problem of acoustic sounding require significant computational resources in the three-dimensional formulation. In this regard, we note articles \cite{12,13}. In them, the initial inverse coefficient problem for the wave equation, from which the equation \eqref{L1_1_} is actually obtained, is reduced to a three-dimensional Fredholm integral equation of the first kind with the right-hand side containing special integrals of the recorded field. The method for solving this equation proposed there turned out to be very effective numerically, and allows solving three-dimensional inverse problems for sufficiently fine grids on a personal computer (PC) even without parallelization. This method is also very useful for the task \eqref{L1_2_}, and this will be demonstrated below.

In this article, we adhere to the following scheme for solving the nonlinear system \eqref{L1_2_}, which reduces the latter to the solution of linear integral equations.

\noindent 1) We solve the second equation written in the form
\begin{equation} \label{L1_3_} \omega ^{2} \int _{X}^{}G(|\textbf{x}-\textbf{x}'|,\omega )v(\textbf{x}',\omega )d\textbf{x}'=W(\textbf{x},\omega ),\, \, \textbf{x}\in Y,
\end{equation}
with respect to the function $v(\textbf{x}',\omega )=\xi (\textbf{x}')u(\textbf{x}',\omega ),\, \, \textbf{x}'\in X$;

\noindent 2) We compute the function $u(\textbf{x},\omega ),\, \, \textbf{x}\in X,$ from the found function $v(\textbf{x}',\omega )$ and from the first equality of the system \eqref{L1_2_}, represented in the form
\begin{equation} \label{L1_4_}
u(\textbf{x},\omega )=u_{0} (\textbf{x},\omega )+\omega ^{2} \int _{X}^{}G(|\textbf{x}-\textbf{x}'|,\omega )v(\textbf{x}',\omega )d\textbf{x}' ,\, \, \, \textbf{x}\in X;
\end{equation}

\noindent 3) We find the function $\xi (\textbf{x})$ from the equation $v(\textbf{x},\omega )=\xi (\textbf{x})u(\textbf{x},\omega ),\, \, \textbf{x}\in X$, using the calculated functions $v(\textbf{x}',\omega )$ and $u(\textbf{x},\omega )$.

A similar scheme was applied, for example, in \cite{Smi}. However, in this paper, regularization methods were not used to solve a three-dimensional equation similar to \eqref{L1_3_}, although this equation is in fact ill-posed. In addition, a solution to the inverse problem was sought on a fairly narrow class of functions with a “piecewise constant current”. Moreover, significant computational resources are required as well to solve such a three-dimensional inverse problem on relatively fine grids.

Below we show that under certain, not very burdensome, special assumptions about the domains $ Y $ and $ X $, the scheme 1) - 3) can be effectively implemented numerically and this allows us to solve the corresponding three-dimensional inverse problem for fairly fine grids in a few minutes even on PCs with medium performance (without parallelization). The proposed algorithm for solving the inverse problem and its numerical study are the main results of this work.

\section{Inverse problem's data recording scheme and task reduction}
We will consider a specific data recording scheme. We assume that the function $u(\textbf{x},\omega )$ is measured “in a flat layer” $Y$ for Cartesian coordinates $\textbf{x}=(x,y,z)$. This scheme was previously used in \cite{12,13} to solve the inverse problem with other data.
\begin{figure}[h]
  \centering
%%%%%%%%%%%%%%%%%%%%%%%%%%%%%%%%%%%%%%%%%%%%%
\includegraphics[width=60mm,height=60mm]{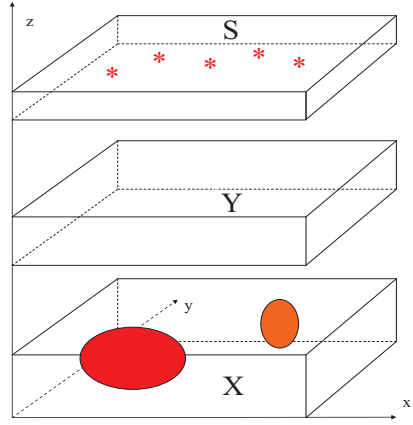}% bmp
%%%%%%%%%%%%%%%%%%%%%%%%%%%%%%%%%%%%%%%%%%%%%
  \caption{{\small Geometric registration scheme of the inverse problem data. Here $X$ is a domain of wave field scatterers, $Y$ is a domain of recording the data, $u(\textbf{x},\omega )$; asterisks show possible positions of the field sources.}}
  \label{fig1}
\end{figure}
Fig.\ref{fig1} shows schematically the corresponding geometric scheme, namely, the region $X$ of the velocity field inhomogeneities, the region $S$ in which the sources are located, and the region of registration of the scattered field $Y$. These areas have the form of infinite flat layers. In more detail, $X={\bf {\mathbb R}}_{xy}^{2} \times [z_{1} ,z_{2} ],\, \, Y={\bf {\mathbb R}}_{xy}^{2} \times [z_{3} ,z_{4} ]$. Asterisks indicate  possible positions of the field sources.

We will use the two-dimensional Fourier transform $F_r[\cdot]$ with respect to the variables $r=(x,y)$, which is defined for the function $a(r,z,\omega )=a(x,y,z,\omega )$ as
\[\tilde{A}(z,\omega ,\Omega )=F_{r} \left[a(r,z,\omega )\right]\left(\Omega \right)=\int _{{\bf {\mathbb R}}_{xy}^{2} }a(x,y,z,\omega ) \exp \left(i(\Omega _{1} x+\Omega _{1} y)\right)dxdy,\, \, \, \, \Omega =(\Omega _{1} ,\Omega _{2} ),\]
as well as the inverse transformation $F_{\Omega }^{-1} \left[\tilde{A}(z,\omega ,\Omega )\right](r)$. Assuming the existence of the corresponding Fourier transforms, i.e. $ \tilde {U}, \, \tilde {U} _ {0}, \, \tilde {G}, \, \tilde {V}, \, \tilde {W} $, for functions $ u, u_ {0}, G, v, W $, we formally converse the integrals
\[\int _{X}^{}G(|\textbf{x}-\textbf{x}'|,\omega )v(\textbf{x}',\omega )d\textbf{x}' ,\, \textbf{x}\in X;\, \, \, \int _{X}^{}G(|\textbf{x}-\textbf{x}'|,\omega )v(\textbf{x}',\omega )d\textbf{x}',\, \, \textbf{x}\in Y ,\]
in the formulas \eqref{L1_3_}, \eqref{L1_4_}. Then we obtain from these formulas the set of equalities
\begin{eqnarray}\label{L1_5_}
&&\tilde U(z,\omega ,\Omega ) = {{\tilde U}_0}(z,\omega ,\Omega ) + {\omega ^2}\int_{{z_1}}^{{z_2}} {\tilde G(z - z',\omega ,\Omega )\tilde V(z',\omega ,\Omega )dz'} ,\,\,\,z \in [{z_1},{z_2}];\nonumber \\
&&\tilde V(z',\omega ,\Omega ) = {F_r}\left[ {\xi (r,z')F_\Omega ^{ - 1}\left[ {\tilde U(z',\omega ,\Omega )} \right](r)} \right]\left( \Omega  \right),\,\,\,z' \in [{z_1},{z_2}]\,;\\
&&{\omega ^2}\int_{{z_1}}^{{z_2}} {\tilde G(z - z',\omega ,\Omega )\tilde V(z',\omega ,\Omega )dz'}  = \tilde W\left( {z,\omega ,\Omega } \right),\,\,\,\,z \in [{z_3},{z_4}].\nonumber
\end{eqnarray}
These equalities are used below to solve the direct and inverse problems. The formulas \eqref{L1_5_} are true if the functions $u,u_{0},G,W$ belong to the class of tempered distribufions \cite{14}, and the function $\xi$ is compactly supported. This will be assumed in what follows. We emphasize that the equalities \eqref{L1_5_} bind the quantities $\tilde U(z,\omega ,\Omega ),\,\tilde G(z - z',\omega ,\Omega )$ and $\tilde V(z',\omega ,\Omega )$ as functions of the arguments $ z, z' $, and the values $\omega ,\Omega$ play the role of parameters.

\section{Algorithm for solving the direct problem}

Our direct problem is formulated as follows. Given a compactly supported function $\xi (r,z')$ and a source function $\tilde{U}_{0} (z,\omega ,\Omega )$, find the value $\tilde{W}\left(z,\omega ,\Omega \right),\, z\in [z_{3} ,z_{4} ]$. To solve this problem, we use the following algorithm.

\textbf{Algorithm 1.}

\noindent 1) Calculate Fourier transforms $\tilde{U}_{0} (z,\omega ,\Omega )$, $\tilde{G}(z-z',\omega ,\Omega )$.

\noindent 2) Implement the iterative process of solving the first two equations from \eqref{L1_5_} for each of the parameters $\omega ,\Omega$ as follows
\begin{multline} \label{L1_6_}
\tilde{U}_{n+1} (z,\omega ,\Omega )=\tilde{U}_{0} (z,\omega ,\Omega )+\omega ^{2} \int _{z_{1} }^{z_{2} }\tilde{G}(z-z',\omega ,\Omega )F_{r} \left[\xi (r,z')F_{\Omega }^{-1}
 \left[\tilde{U}_{n} (z',\omega ,\Omega )\right](r)\right]\left(\Omega
 \right)dz' ,
 \\  z\in [z_{1} ,z_{2} ];\,\,\,n=0,1,2,...,%\nonumber
\end{multline}
starting from the function $\tilde{U}_{0} (z,\omega ,\Omega )$. This can be done in matrix form.

\noindent 3) Stop the iterative process according to some rule at iteration number $\nu $ and get the approximate solution $\tilde{U}_{\nu } (z,\omega ,\Omega )$.

\noindent 4) Sequentially calculate the values
\begin{eqnarray} \label{L1_7_}
&&\tilde{V}_{\nu } (z',\omega ,\Omega )=F_{r} \left[\xi (r,z')F_{\Omega }^{-1} \left[\tilde{U}_{\nu } (z',\omega ,\Omega )\right](r)\right]\left(\Omega \right),\, \, \, z'\in [z_{1} ,z_{2} ] ,
 \nonumber \\ &&\tilde{W}_{\nu } \left(z,\omega ,\Omega \right)={\omega ^2}\int _{z_{1} }^{z_{2} }\tilde{G}(z-z',\omega ,\Omega )\tilde{V}_{\nu } (z',\omega ,\Omega )dz' ,\, \, \, \, z\in [z_{3} ,z_{4} ]
\end{eqnarray}
and take the function $\tilde{W}_{\nu } \left(z,\omega ,\Omega \right)$ or its inverse Fourier transform as an approximate solution to our direct problem.

We will not give here a theoretical analysis of the convergence of Algorithm 1, because formally, we do not use it in solving the inverse problem. This algorithm is only needed to generate model data. We only note that, as follows from the general theory of solving integral equations of the second kind (see, for example, \cite{15}), the algorithm will quickly converge at least for small $\omega $. Numerical examples confirming this statement will be demonstrated below.

\section{Inverse problem}
Our inverse problem is to find the function $\xi (r,z'),\,  z'\in [z_{1} ,z_{2} ]$ from the system \eqref{L1_5_} knowing the value $\tilde{W}\left(z,\omega ,\Omega \right),\,  z\in [z_{3} ,z_{4} ]$. To do this, we apply the following algorithm.

\textbf{Algorithm 2.}

\noindent 1) Given a function $\tilde{W}\left(z,\omega ,\Omega \right)$, we solve one-dimensional Fredholm equations of the first kind
\begin{equation}\label{L1_8_}
\left\{{\omega ^2}\!\int _{z_{1} }^{z_{2} }\tilde{G} (z-z',\omega ,\Omega )\tilde{V}(z',\omega ,\Omega )dz' =\tilde{W}\left(z,\omega ,\Omega \right),\, \, \, \, z\in [z_{3} ,z_{4} ]\right\}\, \, \, \Rightarrow \tilde{V}(z',\omega ,\Omega ),
\end{equation}
for all considered parameters $\omega ,\Omega $. To do this, we use the appropriate regularization method (regularizing algorithm, RA) for these ill-posed problems.

\noindent 2) Using the found function $\tilde{V}(z',\omega ,\Omega )$, we calculate the quantity
\begin{equation}\label{L1_9_}
\tilde{U}(z,\omega ,\Omega )=\tilde{U}_{0} (z,\omega ,\Omega )+\omega ^{2} \int _{z_{1} }^{z_{2} }\tilde{G}(z-z',\omega ,\Omega )\tilde{V}(z',\omega ,\Omega )dz' ,\, \, \, z\in [z_{1} ,z_{2} ],
\end{equation}
and further the functions
\[V(r,z',\omega )=F_{\Omega }^{-1} \left[\tilde{V} (z',\omega ,\Omega )\right](r),\, \, \, \, u(r,z',\omega )=F_{\Omega }^{-1} \left[\tilde{U}(z',\omega ,\Omega )\right](r),\, \, \, (r,z')\in X\, .\]
\noindent 3) We find the solution $\xi (r,z')$ from the equation $u(r,z',\omega )\xi (r,z')=V(r,z',\omega )$. This can be done using the least squares method with respect to $\omega $ or simply calculating the solution in the form $\xi (r,z')=V(r,z',\omega )/u(r,z',\omega )$. In the latter case, the result will, generally speaking, depend on $\omega $.

Let's make some comments about Algorithm 2.

a) To implement item 1, it is necessary to refine assumptions about functions ${V},\,{W}$. Namely, we use inclusions ${V}(\textbf{x},\omega)\in L_2(\mathbb{R}^3_{xyz}),\,{W}(\textbf{x},\omega)\in L_2(\mathbb{R}^3_{xyz})$ for each considered frequency $\omega$. The first inclusion follows from the finiteness of the function $\xi$, and the second is postulated. In this case, many methods for solving linear ill-posed problems are applicable to the equations \eqref{L1_8_} (see, for example, \cite{2,3,16,17,18} and others).

b) The equation \eqref{L1_3_} and consequently, the equations \eqref{L1_8_}, may have more than one solution for the finite set of frequencies used. Therefore, to implement the algorithm, we use RAs that give normal solutions to the equations (for example, Tikhonov regularization, TSVD method). If the solution to the equation \eqref{L1_3_} is unique, it will coincide with the calculated normal solution.

c) The equations \eqref{L1_8_} can be solved separately for each fixed frequency $\omega$, or they can also be solved as a system of equations for all frequencies used. Accordingly, in the latter case, item 1 of the algorithm is implemented using least squares.

\section{Numerical experiments}
In this section, it is assumed that the equations \eqref{L1_5_} are written in dimensionless form with $c_{0} =1$, so $k_{0} =\omega $. We also assume that model sources are given in the form $f(\textbf{x})=\sum _{m}A_{m} \delta (\textbf{x}-\textbf{x}_{m} ) $, where $\textbf{x}_{m} $ are the coordinates of $\delta $-shaped sources. Then $u_{0} (\textbf{x},\omega )=\sum _{m}A_{m} G(|\textbf{x}-\textbf{x}'_{m} |,\omega ) $, and the Fourier transform of this function, $\tilde{U}_{0} (z,\omega ,\Omega )$, can be calculated.

For a numerical study of the proposed algorithms, we use following model solutions:
\[\xi (x,y,z)=A_{0} \left[\left(1-\frac{\rho _{1} ^{2} }{0.4^{2} } \right)_{+} +2\left(1-\frac{\rho _{2} ^{2} }{0.25^{2} } \right)_{+} +2.5\left(1-\frac{\rho _{3} ^{2} }{0.3^{2} } \right)_{+} \right]\]
Here \vspace{-2mm}
\begin{eqnarray*}
  &&\rho _{1} =\sqrt{(x-1)^{2} +(y-2)^{2} +(z-0.5)^{2} } ,\\
  &&\rho _{2} =\sqrt{(x-4)^{2} +(y+3)^{2} +(z-0.5)^{2} +1.5(y+3)(z-0.5)} ,\\
  &&\rho _{3} =\sqrt{(x+3)^{2} +y^{2} +(z-0.45)^{2} -1.5y(z-0.45)} , \,\, \left(a\right)_{+} =\max \left(a,0\right).
\end{eqnarray*}
These functions correspond to small local inhomogeneities of the medium, the position of which and the corresponding velocity distributions must be found. It is to search for such heterogeneities that Algorithm 2 is configured.

The value $A_{0} $ determines the so-called contrast $\frac{\Delta c}{c_{0} } =\max \left\{\frac{1}{\sqrt{1-c_{0}^{2} \xi (r)} } \right\}-1$ of the  solution to be found. In the calculations presented below, the value $A_{0} =0.3$ was used, which corresponds to a contrast of 25.2. According to the classification from \cite[c.33]{1}, such scatterers can be considered as “strong” because $\frac{\Delta c}{c_{0} }\gg \frac{c_0}{l\omega}$ for their characteristic dimensions $l\sim 0.4$ and the values $c_0=1,\omega=2$. Scatterers of this type are quite common in practice.

The initial equations \eqref{L1_3_}, \eqref{L1_4_} were approximated by the finite-difference method on uniform grids in domains
\[X=[-10,10]\times [-10,10]\times [-0.5,1.5],\,Y=[-10,10]\times [-10,10]\times [6.01,6+\varepsilon ].\]
The dimensions of the grids in these areas were chosen as $N\times N\times M$ and $N\times N\times M_{1} $ with the numbers $N,M,M_{1} $, which will be indicated below. The number $\varepsilon >0$ characterizes the thickness of the layer in which the data is measured. We consider two variants,  $\varepsilon =0.5$ (thick layer) or $\varepsilon =0.02$ (thin layer).

Discrete analogues of the functions $\tilde{U}_{0} (z,\omega ,\Omega )$, $\tilde{G}(z-z',\omega ,\Omega )$ were calculated for the given frequencies $\omega $ from the known values of $u_{0} (x,\omega )$ and $G(\rho ,\omega )$ using the fast Fourier transform (FFT). Then they were used in Algorithms 1 and 2.

All model problems were solved with various sources. In this paper, we do not set as our goal the optimization of their position and number. We only note that for the solved problems with sources located at points \[\textbf{x}_{m} =(0,y_{m} ,6),\, y_{m} =[-5,-4,...,0,...,4,5]\] with $ A_{m} = 1 $ and with similar sources at points \[\textbf{x}_{m} =\, \, ({\rm 0,0,6),\; (}\pm {\rm 5,0,6),\; }(0{\rm ,}\pm {\rm 5,6),\; (}\pm {\rm 2,0,6),\; }(0{\rm ,}\pm {\rm 2,6),\; (}\pm 1{\rm ,0,6),\; }(0{\rm ,}\pm 1{\rm ,6)},\] the results were very close. Therefore, we present below the results of numerical experiments only for sources of the first kind.

\subsection{Using Algorithm 1 to obtain model data for solving the inverse problem}
Here we present typical results of a numerical study of the iterative process \eqref{L1_6_} to obtain the data of the inverse problem. The convergence rates of the process \eqref{L1_6_} for various quantities $\omega =k_{0} $ are compared in Fig.2. Here we use the notation $u_{n} =\tilde{U}_{n} (z,\omega ,\Omega )$.
\begin{figure}[h]
  \centering
%%%%%%%%%%%%%%%%%%%%%%%%%%%%%%%%%%%%%%%%%%%%%
\includegraphics[width=80mm,height=60mm]{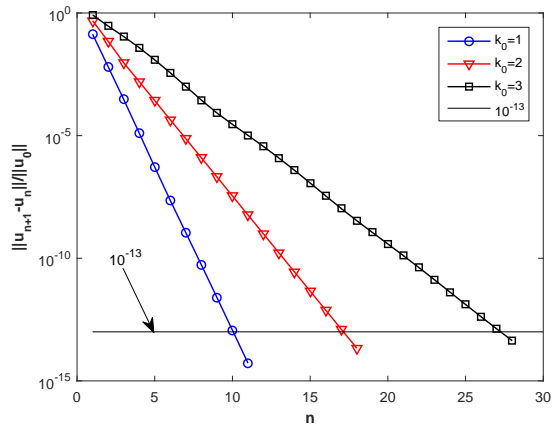}% bmp
%%%%%%%%%%%%%%%%%%%%%%%%%%%%%%%%%%%%%%%%%%%%%
  \caption{{\small The convergence rate of iterations \eqref{L1_6_} for different $\omega =k_{0} $.}}
  \label{fig2}
\end{figure}
\noindent Iterations were stopped when the condition $\left\| \tilde{U}_{\nu } (z,\omega ,\Omega )-\tilde{U}_{\nu -1} (z,\omega ,\Omega )\right\| \le 10^{-13} \left\| \tilde{U}_{0} (z,\omega ,\Omega )\right\| $ holds true. After that, using the found value $\tilde{U}_{\nu } (z,\omega ,\Omega )$, we calculated by the formulas \eqref{L1_7_} the function $\tilde{W}_{\nu } \left(z,\omega ,\Omega \right)$ and then the function \[W_{\nu } (x,y,z,\omega )=F_{\Omega }^{-1} \left[\tilde{W}_{\nu } \left(z,\omega ,\Omega \right)\right](r), \, z\in [z_{3} ,z_{4} ].\] The last function provides data for solving the inverse problem. Its form, found for $\omega =k_{0} =2$, is shown in Fig.3 for $ z = 6.25 $.
\begin{figure}[h]
  \centering
%%%%%%%%%%%%%%%%%%%%%%%%%%%%%%%%%%%%%%%%%%%%%
\includegraphics[width=100mm,height=50mm]{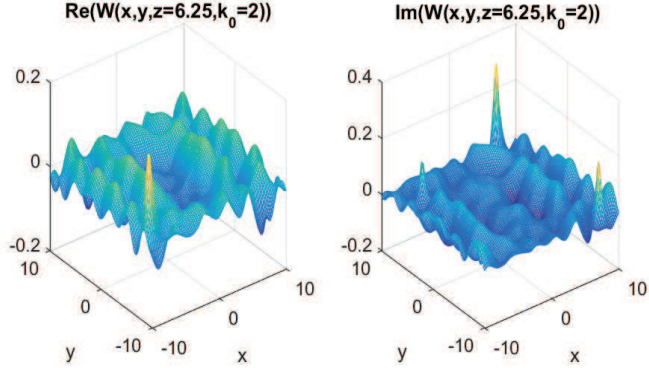}% bmp
%%%%%%%%%%%%%%%%%%%%%%%%%%%%%%%%%%%%%%%%%%%%%
  \caption{{\small Typical data for solving the inverse problem.}}
  \label{fig3}
\end{figure}
Generally speaking, the data of the inverse problem can be given with some perturbations, which are interpreted as measurement errors. In our calculations, this was modeled by superimposing on the function $W_{\nu } (x,y,z,\omega )$ an additive normally distributed pseudo-random noise with a zero mean so that the resulting approximate function $W_{\nu }^{(\delta )} (x,y,z,\omega )$ would satisfy the condition
\[\left\| W_{\nu }^{(\delta )} (x,y,z,\omega )-W_{\nu } (x,y,z,\omega )\right\| _{L_{2} ({\bf {\mathbb R}}_{xyz}^{3} )} \le \delta \left\| W_{\nu } (x,y,z,\omega )\right\| _{L_{2} ({\bf {\mathbb R}}_{xyz}^{3} )}. \] This corresponds to approximate data with relative accuracy $\delta $.

\subsection{Implementation of Algorithm 2 }
Item 1 of Algorithm 2, i.e. the solution of equations of the first kind \eqref{L1_8_} by regularization methods, was discussed in \cite{12,13} in connection with the solution of another inverse problem. It was formally similar to the problem under consideration and differing from it only in the form of the kernel and the right-hand side. In these works, it was noted that the discretization used reduces the equations \eqref{L1_8_} to a set of systems of linear algebraic equations (SLAE) of the form $A^{(m)} \tilde{V}^{(m)} =\tilde{W}^{(m)} $ with matrices $A^{(m)} =\left[\mu _{kl} \tilde{G} (z_{k} -z_{l} ^{{'} } ,\omega ,\Omega ^{(m)} )\right]$ of size $M\times M_{1} $ and the right parts $\tilde{W}^{(m)} =\left[\tilde{W}_{\nu }^{(\delta)} \left(z_{k} ,\omega ,\Omega ^{(m)} \right)\right]$, which are column vectors of height $ M $. Here $\Omega ^{(m)} $ are nodes of a two-dimensional grid of $\Omega $ quantities numbered by a single index $ m $, and $\mu _{kl} $ are quadrature coefficients for calculating the integrals in \eqref{L1_8_}. We solved the mentioned SLAEs using various variants of the Tikhonov regularization \cite{16,17} and using the TSVD method \cite{18}. The best results in the calculations were obtained for the TSVD method. We present them below.

Item 2 of Algorithm 2 does not cause difficulties for a discretized problem, because reduces to matrix multiplication of discrete quantities $\tilde{G},\tilde{V}$ and addition of the result with a discrete analogue of the function $\tilde{U}_{0} $. Further calculations of this item are performed using the FFT.

Finally, item 3 is realized by dividing $\xi (r,z')=V(r,z',\omega )/u(r,z',\omega )$ in the case of solving the inverse problem for a single frequency $\omega $ or its separate solution for several frequencies. In the case of a joint solution for several frequencies, we use least squares in the variable $\omega $. In the end, if necessary, it is easy to recalculate the function $\xi(\mathbf{x})$ into $c(\mathbf{x})$. For the sake of brevity, we do not do this in the examples below, presenting the value $\xi(\mathbf{x})$ in the figures.
\begin{figure}%[h]
  \centering
%%%%%%%%%%%%%%%%%%%%%%%%%%%%%%%%%%%%%%%%%%%%%
\includegraphics[width=160mm,height=80mm]{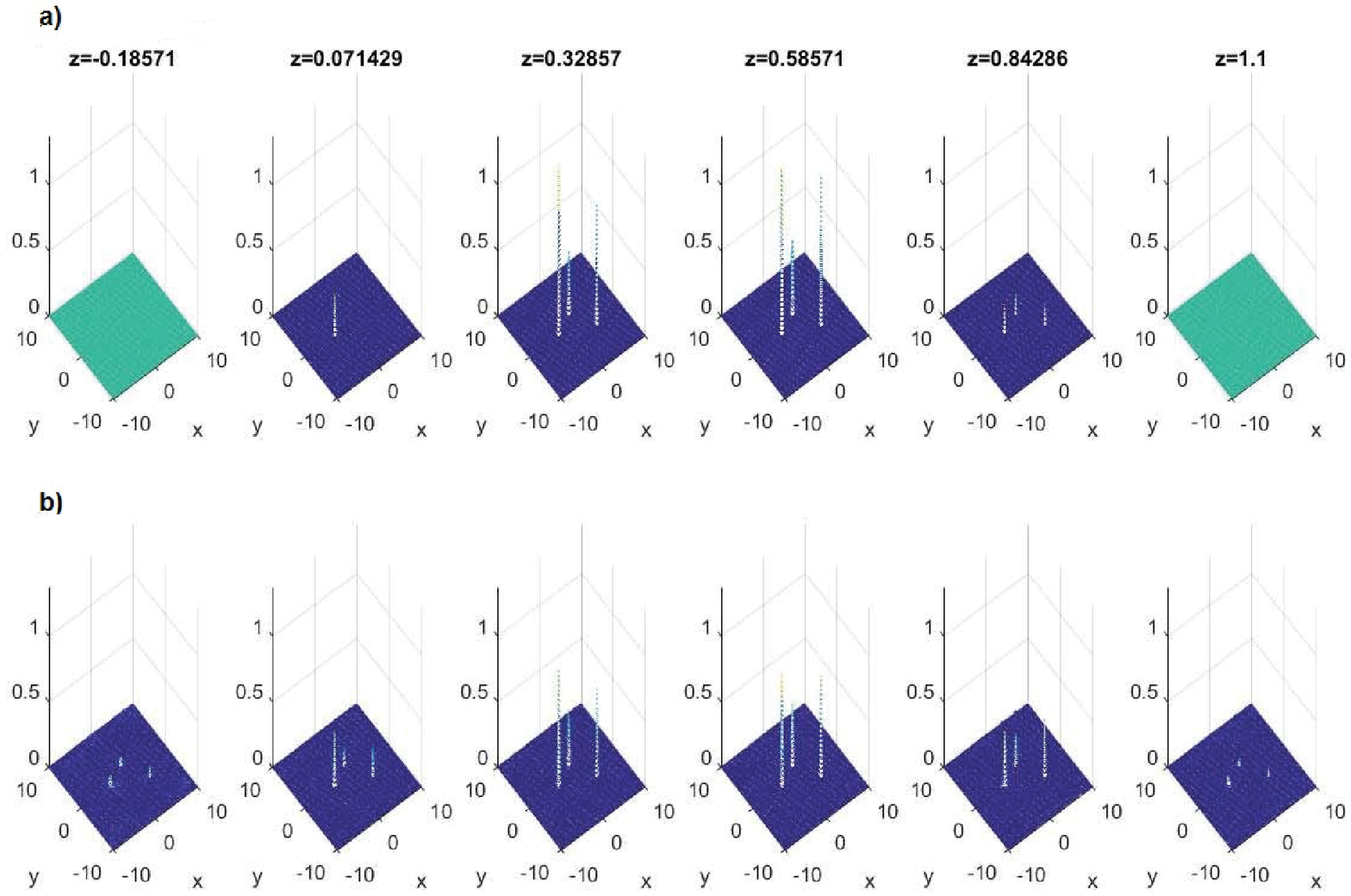}% bmp

\includegraphics[width=160mm,height=40mm]{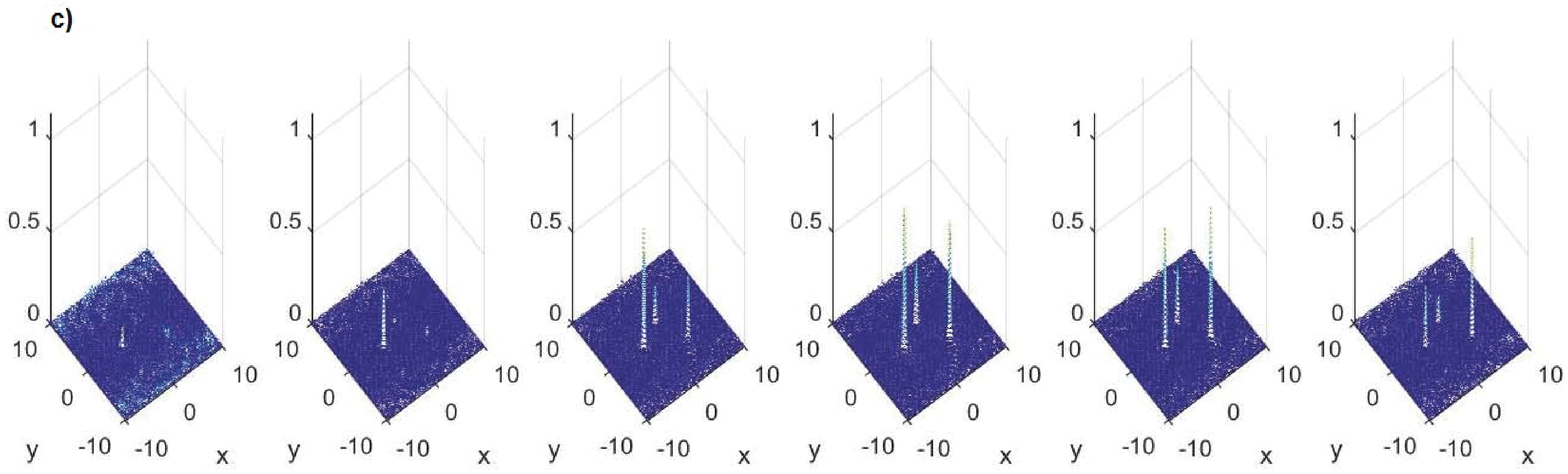}
\includegraphics[width=160mm,height=40mm]{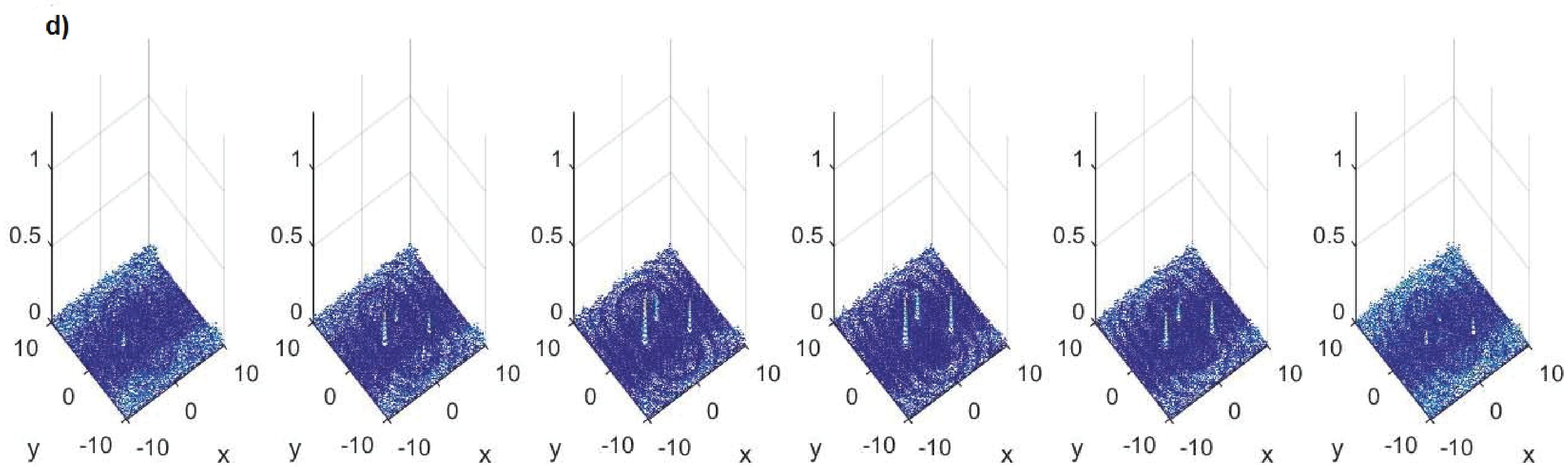}
%%%%%%%%%%%%%%%%%%%%%%%%%%%%%%%%%%%%%%%%%%%%%
  \caption{{\small A qualitative comparison of the exact (a) and approximate solutions of the inverse problem: b) for exact data (in the thick layer); c) for perturbed data with $\delta =10^{-7} $; d) for perturbed data with $\delta =10^{-5} $.}}
  \label{fig4}
\end{figure}

\subsection{Solving the inverse problem with data in a thick layer.}

In this group of experiments, the data of the inverse problem are specified in the layer $Y=[-10,10]\times[-10,10]\times[6.01,6.5]$ at one frequency, $\omega=k_0=2$. The grids of sizes $N=128,M=M_1=71$ are used. Fig.4 makes it possible to compare qualitatively the exact and approximate solutions, $\xi _{exact} (x,y,z)$ and $\xi _{appr} (x,y,z)$, of the inverse problem for various $ z $. The exact solution is shown in the first line of the figure. The second line gives an approximate solution of the inverse problem for exact data. The third and fourth lines depict solutions for perturbed data with $\delta =10^{-7} $ and $\delta =10^{-5} $, respectively. The figure demonstrates a fairly high sensitivity of the solutions to data perturbations. %Solutions for more noisy data are smoother.

More detailed information on the accuracy of solving the inverse problem, i.e. about relative accuracy
\[\Delta _{L_{2} } (z)=\frac{\left\| \xi _{appr} (x,y,z)-\xi _{exact} (x,y,z)\right\| _{L_{2} ({\bf {\mathbb R}}_{xy}^{2} )} }{\left\| \xi _{exact} (x,y,z)\right\| _{L_{2} ({\bf {\mathbb R}}_{xy}^{2} )} } \]
of the approximate solution $\xi _{appr} (x,y,z)$ for different $z$ and $\delta $ are presented in Fig.5.
\begin{figure}%[h]
  \centering
%%%%%%%%%%%%%%%%%%%%%%%%%%%%%%%%%%%%%%%%%%%%%
\includegraphics[width=80mm,height=70mm]{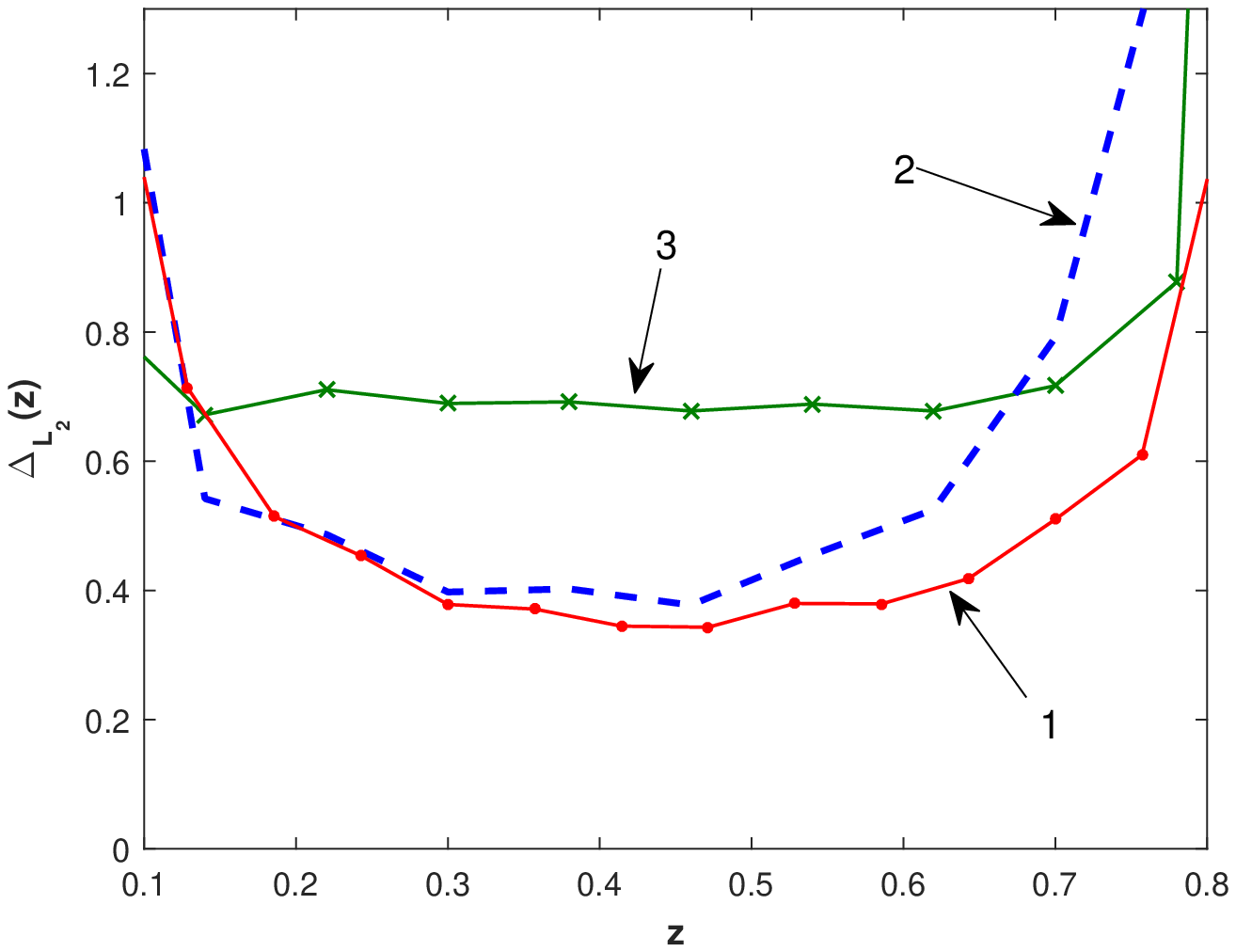}% bmp
%%%%%%%%%%%%%%%%%%%%%%%%%%%%%%%%%%%%%%%%%%%%%
  \caption{{\small Relative accuracy $\Delta _{L_{2} } (z)$ of an approximate solution (in the norm of the space $L_{2} ({\bf {\mathbb R}}_{xy}^{2} )$ for various $ z $. Line 1 -- value $\Delta _{L_{2} } (z)$  for unperturbed data of the problem, line 2 -- for approximate data with $\delta =10^{-7} $, line 3 -- for approximate data with $\delta =10^{-5} $.}}
    \label{fig5}
\end{figure}

Fig.6 shows qualitatively the influence of a perturbation of the problem data on determining the position of inhomogeneities.
\begin{figure}%[h]
  \centering
%%%%%%%%%%%%%%%%%%%%%%%%%%%%%%%%%%%%%%%%%%%%%
\includegraphics[width=130mm,height=80mm]{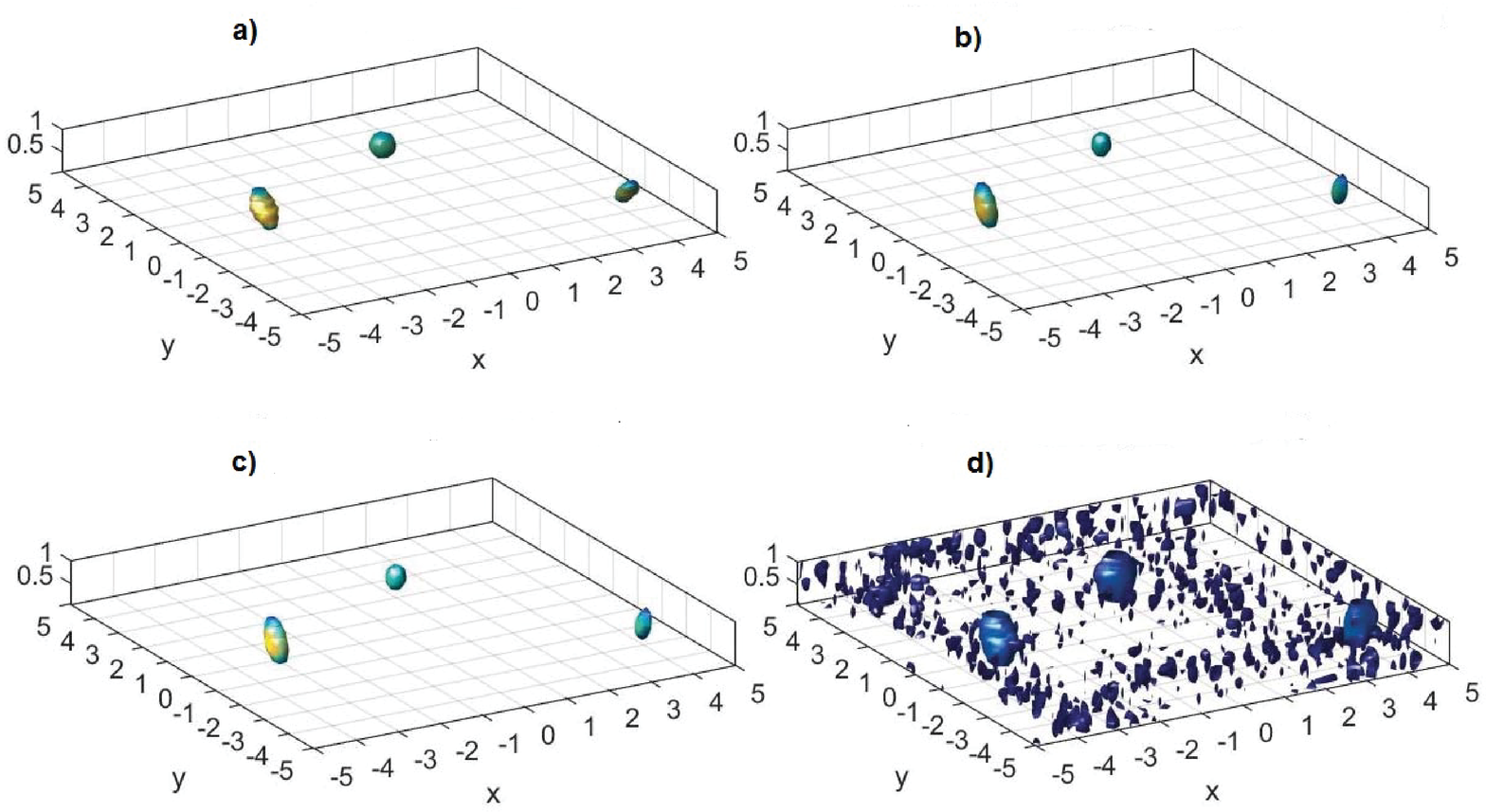}% bmp
%%%%%%%%%%%%%%%%%%%%%%%%%%%%%%%%%%%%%%%%%%%%%
  \caption{{\small Qualitative comparison of the influence of data perturbation on the position of the reconstructed inhomogeneity $\xi(\mathbf{x})$ ($k_{0} =2$). a) exact solution, b) approximate solution for exact data, c) approximate solution for $\delta=10^{-7}$, d) approximate solution for $\delta=10^{-5}$}}
  \label{fig6}
\end{figure}
It is seen that the positions can be quite accurately determined in the case of $\delta =10^{-7} $, and for $\delta =10^{-5} $ this is possible using appropriate noise filtering.

Now we demonstrate the solutions of the inverse problem for different frequencies $\omega=k_0=1,2,3$. In this case, we use the exact data of the problem in order to more clearly highlight how the frequency affects the results. Fig.7 shows how the positions of inhomogeneities are restored when solving the inverse problem separately for each frequency from the set.
\begin{figure}%[h]
  \centering
%%%%%%%%%%%%%%%%%%%%%%%%%%%%%%%%%%%%%%%%%%%%%
\includegraphics[width=130mm,height=80mm]{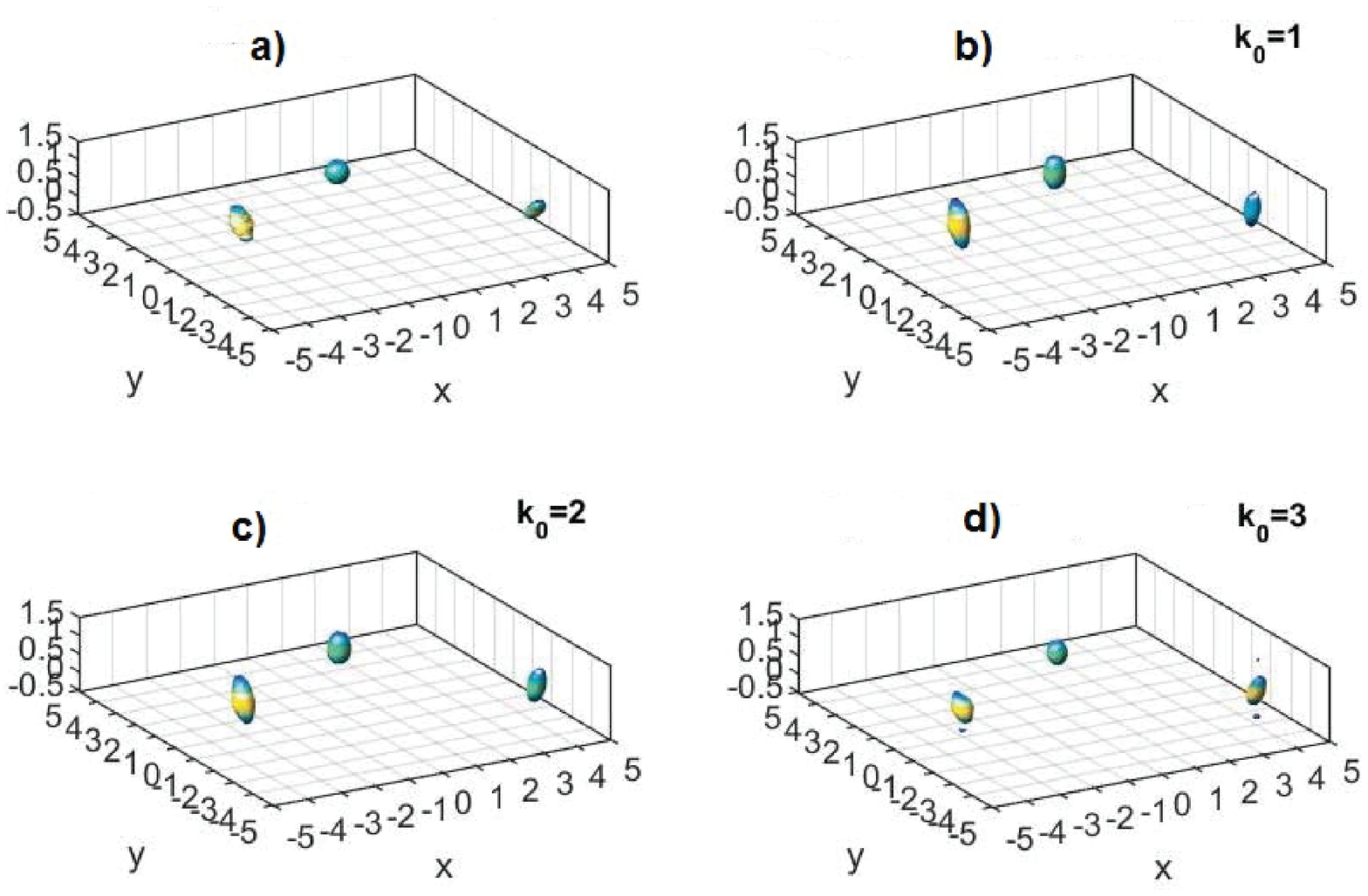}% bmp
%%%%%%%%%%%%%%%%%%%%%%%%%%%%%%%%%%%%%%%%%%%%%
  \caption{{\small A qualitative comparison of the positions and geometry of the reconstructed inhomogeneity $\xi(\mathbf{x})$ for various $k_0$ (unperturbed data); a) exact solution; b) approximate solution for $k_0=1$, c) for $k_0=2$, d) for $k_0=3$.}}
   \label{fig7}
\end{figure}

The inverse problem can also be solved by considering the equalities \eqref{L1_8_} and \eqref{L1_9_} as a system of equations for several frequencies. Figure 8 presents the results of solving the inverse problem for three frequencies at once: $\omega =k_{0} =[1,2,3]$, with finding the final function $\xi (r,z')$ from the equation $u(r,z',\omega )\xi (r,z')=V(r,z',\omega )$ according to the least-squares method. %Note that the result turned out to be slightly worse than, e.g., for $k_0=2$.
\begin{figure}%[h]
  \centering
%%%%%%%%%%%%%%%%%%%%%%%%%%%%%%%%%%%%%%%%%%%%%
\includegraphics[width=160mm,height=90mm]{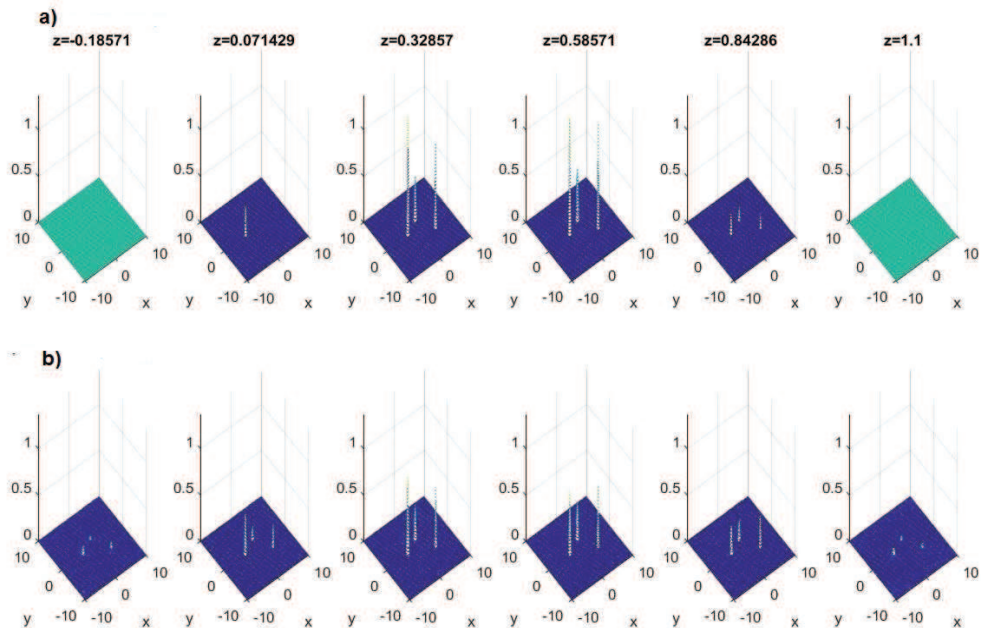}% bmp
%%%%%%%%%%%%%%%%%%%%%%%%%%%%%%%%%%%%%%%%%%%%%
  \caption{{\small Qualitative comparison of the exact solution $ \xi (\mathbf{x}) $ (a) and the approximate solution of the inverse problem simultaneously at three frequencies (b) (for the exact data).}}
  \label{fig8}
\end{figure}
These results are shown in more detail in Fig.9 on the left, where the accuracies of the obtained solutions are compared separately for each of the frequencies and for simultaneous solution of the inverse problem at three frequencies. One can see that the least squares result (green line) is slightly worse than, e.g., for $k_0=2$ (red line). This can be explained. Algorithm 2, when applied for each frequency, gives different accuracies of approximate solutions. Therefore, the accuracy of the joint solution at several frequencies with obtaining the final result by the least squares method may turn out to be worse than at some specific frequencies from the set used. It is for this reason that the inverse problem should be solved separately for each frequency, and not just for their combination at once.
%For example, one can see in Fig.9 that the accuracy of the solution for $k_0=2$ is generally better than the accuracy of the joint solution.
\begin{figure}%[h]
  \centering
%%%%%%%%%%%%%%%%%%%%%%%%%%%%%%%%%%%%%%%%%%%%%
\includegraphics[width=80mm,height=60mm]{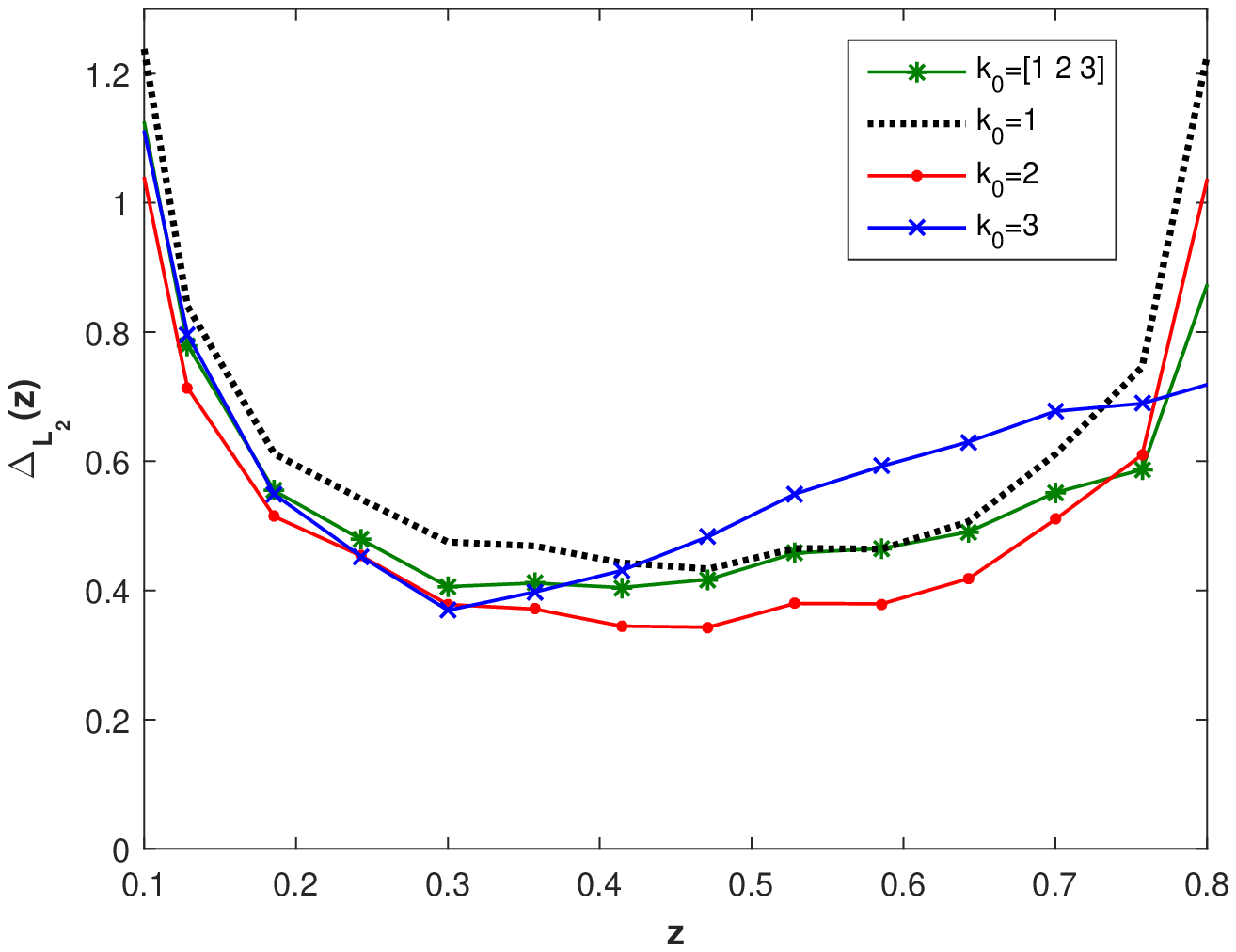}% bmp
\includegraphics[width=80mm,height=60mm]{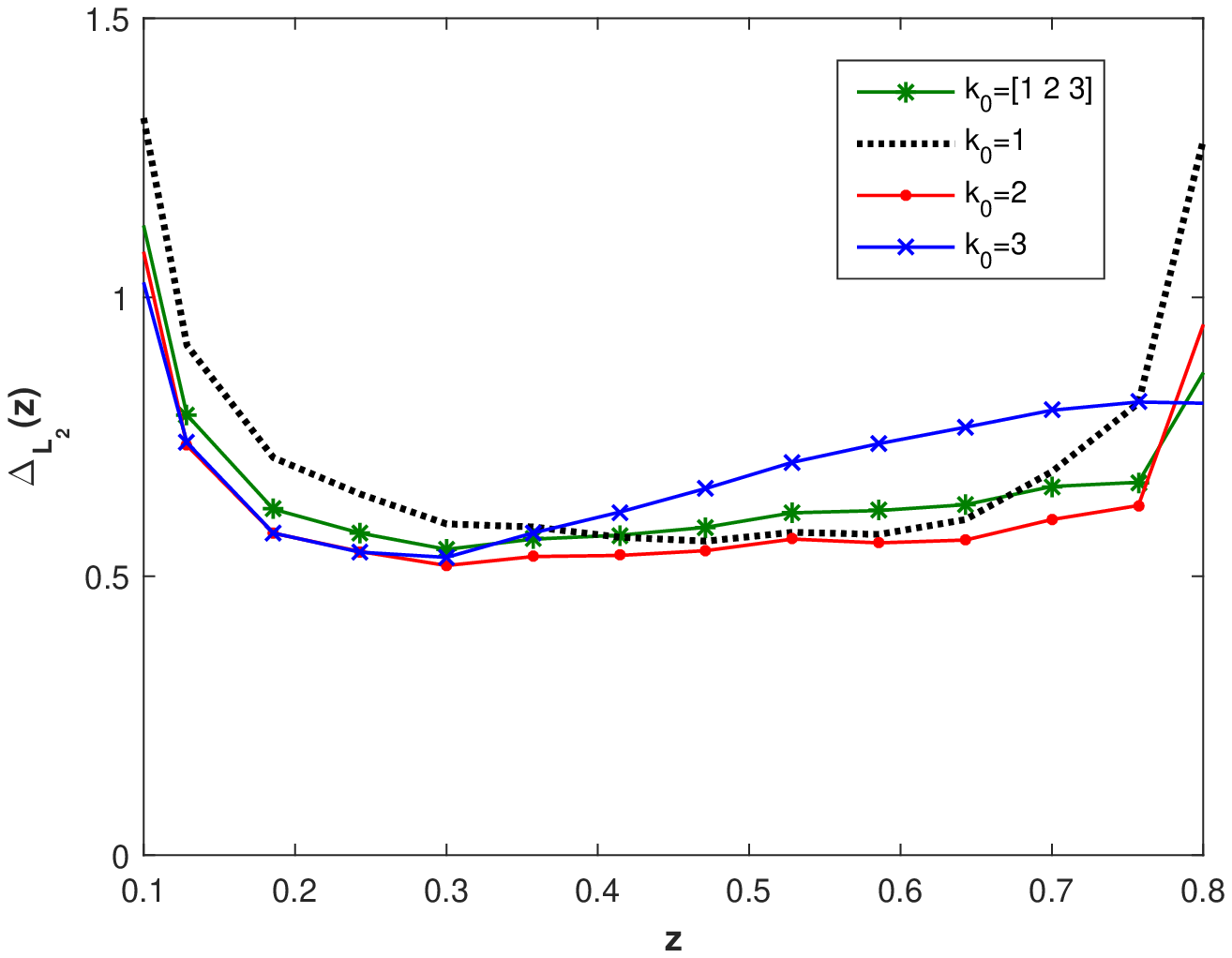}
%%%%%%%%%%%%%%%%%%%%%%%%%%%%%%%%%%%%%%%%%%%%%
  \caption{{\small Left: relative accuracy of approximate solutions, $\Delta _{L_{2} } (z)$, with data in a thick layer for different $k_{0} $, and simultaneously for all frequencies $k_{0} =[1,2,3]$ (unperturbed data). Right: similar accuracy of approximate solutions obtained from unperturbed data in a thin layer.}}
    \label{fig9}
\end{figure}

\subsection{Solving the inverse problem for data in a thin layer}
In these experiments, uniform grids were used in the regions
\[X=[-10,10]\times [-10,10]\times [-0.5,1.5],\, \, \, Y=[-10,10]\times [-10,10]\times [6.01,6.02].\] The sizes of the grids were as follows:
 $ N = 128 $ for the variables $ x$ and $y $ in the domains $ X,Y $, $ M = 71 $ for the variable $ z \in [-0.5,1.5] $,  and $ M_1 = 2 $ for the variable $ z'\in [6.01,6.02] $.

The results are summarized in Fig.9 on the right. The calculations were carried out separately for different $k_{0} $ and simultaneously at all frequencies used (green line) using unperturbed data of the inverse problem. It can be seen that the solution error is somewhat worse than when solving in a thick layer. However, even for such data, the algorithm makes it possible to determine reliably the positions of local inhomogeneities.

\section{Some properties of Algorithm 2}
All calculations were carried out in MATLAB on a PC with a processor Intel (R) Core (TM) i7-7700 CPU 3.60 GHz, 16GB RAM without parallelization.
Algorithm 2 for solving the inverse problem turned out to be quite fast. In this connection, we present the results of numerical experiments for solving the inverse problem at one and the same frequency $\omega=2 $ for different grids. The numbers $ M, M_1 $ determine the speed of solving the one-dimensional integral equation \eqref{L1_8_} for a fixed $ \Omega $. 
The corresponding solution time, $t_{0} (M,M_1)$, varies a little in
passing from one of equations to another. This time is controlled by the desired resolution of the algorithm
in the variable $ z $. Actually, we have the estimate $t(N,M,M_1)\approx t_{0} (M,M_1)\cdot N^{2} $ for full solution time of the
inverse problem for chosen grids, and the number $ N $ is controlled here by the required resolution in $ x, y $.
%The corresponding solution time, $t_{0} (M,M_1)$, practically does not change much 
%for various equations of this type. This time is controlled by the required 
%resolution of the algorithm with respect to the variable $ z $. Thus, the total 
%time to solve the inverse problem on the indicated grids can actually be estimated 
%as $t(N,M,M_1)\approx t_{0} (M,M_1)\cdot N^{2} $, and the number $ N $ is controlled 
%here by the required resolution along the axes $ x, y $.

We present in Fig.10 this dependence $t(N)=t(N,M,M_1)$, obtained in experiments for fixed $M=51,M_1=51$. This time can be significantly reduced by using parallel computations.
% и для различных $N$.
\begin{figure}[h]
  \centering
%%%%%%%%%%%%%%%%%%%%%%%%%%%%%%%%%%%%%%%%%%%%%
\includegraphics[width=70mm,height=60mm]{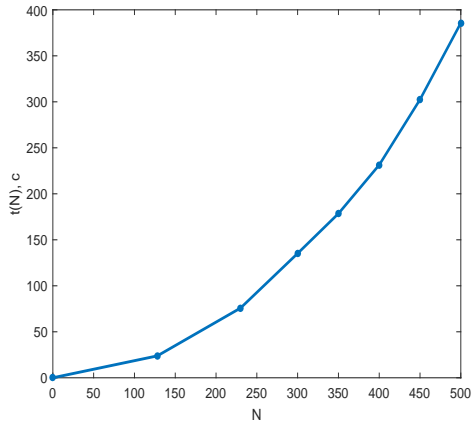}% bmp
%%%%%%%%%%%%%%%%%%%%%%%%%%%%%%%%%%%%%%%%%%%%%
  \caption{{\small The solving time vs $N$.}}
  \label{fig11}
\end{figure}

We also note that the inverse problem under consideration is very sensitive to input data errors. When solving it with double precision, random errors with an amplitude of the order of $10^{-5} $ in the right-hand side of the equation \eqref{L1_8_} lead to significant distortions in approximate solutions obtained by various regularization methods. This is due to the very rapid decay of the singular numbers of the matrices $A^{(m)} $ in SLAEs solved in Algorithm 2, and this is a specific feature of the inverse problem being solved. A similar property of the inverse coefficient problem for the wave equation was noted earlier in \cite{12,13}. Corresponding theoretical estimates of the solution error under various a priori assumptions on the exact solution can be found in \cite{2,3}.

\section{Conclusions}
From the numerical experiments carried out in the work, we can draw the following conclusions.

1. Three-dimensional inverse coefficient problem for the wave equation, which arises in the modeling of the acoustic sensing, can be solved numerically using the proposed Algorithm 2 for sufficiently fine mesh in a few minutes on a typical PC, even without parallelization. To achieve this, one can use, for example, the input data of the inverse problem in a flat layer.

2. By itself, the inverse problem under consideration is very sensitive to data perturbations. To obtain a detailed approximate solution, data measured with great accuracy are required. This feature of the problem does not depend on the algorithm used to solve it.

3. The proposed algorithm makes it possible to determine reliably the positions of small local inhomogeneities of the acoustic medium for data with small errors. Data can be specified even in a thin layer.

4. Algorithm 2, when applied separately for different frequencies, gives different accuracy of approximate solutions. Therefore, the accuracy of the joint solution of the inverse problem at several frequencies according to the least-squares method may turn out to be somewhat worse than at some specific frequencies from the set used.

\bigskip
\textbf{Funding:} The work of the second author was supported by the Programm of Competitiveness Increase of the National Research Nuclear University
MEPhI (Moscow Engineering Physics Institute); contract no. 02.a03.21.0005, 27.08.2013.


\begin{thebibliography}{2}
\bibitem{1} Goryunov A.A., Saskovets A.V. Inverse scattering problems in acoustics. M., Publishing House of Moscow State University, 1980.

\bibitem{21} D. Colton and R. Kress, Inverse Acoustic and Electromagnetic Scattering Theory, 2nd ed., Appl. Math. Sci. 93, Springer,
Berlin, 1998.

\bibitem{22} A. G. Ramm, Multidimensional Inverse Scattering Problems, Pitman Monogr. Surv. Pure Appl. Math. 51, Longman Scientific \& Technical, Harlow, 1992.

\bibitem{2} Bakushinsky А., Goncharsky А. Ill-Posed Problems: Theory and Applications. Dordrecht, Kluwer Academic Publishers, 1994.

\bibitem{3} Bakushinsky A.B., Kokurin M.Yu. Iterative methods for approximate solution of inverse problems. Mathematics and Its Applications. Dordrecht, Kluwer Academic Publishers, 2004.

\bibitem{4} A. V. Goncharsky and S. Y. Romanov, On two approaches to the solution of coefficient inverse problems for wave
equations, Zh. Vychisl. Mat. Mat. Fiz. 52 (2012), no.2, pp.263–269.

\bibitem{5} A. V. Goncharsky and S. Y. Romanov, Supercomputer technologies in inverse problems of ultrasound tomography,
Inverse Problems 29 (2013), no.7, Article ID 075004.

\bibitem{6} Belishev M.I. Recent progress in the boundary control method, Inverse Problems. Vol.23, no.5. 2007, pp.1--67.

\bibitem{7} Pestov L.N., Bolgova V.M., Danilin A.N. Numerical reconstruction of the three-dimensional speed of sound by the method of boundary control, Bulletin of Ugra State University, 2011, Issue 3, pp.92--98 (in Russian).
    
\bibitem{9} Novikov P. G. Reconstruction of the two-dimensional Schrodinger operator from the scattering amplitude at a fixed energy, Funktsional. analysis and its adj., T.20, No.3, 1986, pp.90-91.    

\bibitem{8} Burov V.A.,  Alekseenko N.V.,  Rumyantseva O.D. Multifrequency Generalization of the Novikov Algorithm for the Two-Dimensional Inverse Scattering Problem // Acoustic Journal, V.55, No.6, 2009, pp.84-798.

\bibitem{10} Burov V.A., Vecherin S.N., Morozov S.A., Rumyantseva O.D. Modeling of the Exact Solution of the Inverse Scattering Problem by Functional Methods. Acoustic Journal, Vol.56, No.4, 2010, pp.516-536.

\bibitem{11} Beilina L. and Klibanov M.V., Approximate Global Convergence and
Adaptivity for Coefficient Inverse Problems. New York, Springer, 2012. %408 p.

\bibitem{Kab} Kabanikhin S.I., Satybaev A.D., Shishlenin M.A. Direct Methods of Solving Multidimensional Inverse Hyperbolic Problems. Utrecht, VSP, 2004.

\bibitem{Kl1} Michael V. Klibanov and Aleksandr E. Kolesov. Convexification of a 3-D coefficient inverse scattering
problem, Computers and Mathematics with Applications 77 (2019), pp. 1681–1702

\bibitem{Kl2} Michael V. Klibanov, Aleksandr E. Kolesov, and Dinh-Liem Nguyen. Convexification method for an inverse scattering problem and its performance for experimental backscatter data for buried targets. SIAM J. Imaging Sciences, Vol. 12, No. 1 (2019), pp. 576--603

\bibitem{12} Bakushinsky A.B., Leonov A.S. Fast numerical method of solving 3D coefficient inverse problem for wave equation with integral data, Journal of Inverse and Ill-Posed Problems, 2018, Vol. 26, Issue 4. pp. 477-492.

\bibitem{13} A.B. Bakushinskii, and A.S. Leonov, Low-Cost Numerical Method for Solving a Coefficient Inverse Problem for the Wave Equation in Three-Dimensional Space. Comp. Math. and Math. Phys., 2018, Vol. 58, No. 4, pp.548–561.

\bibitem{Smi}  Evstigneev R.O., Medvedik M.Yu., Smirnov Yu.G., Tsupak A.A. The inverse problem of body’s heterogeneity recovery for early diagnostics of diseases using microwave tomography. University proceedings. Volga region. Physical and Mathematical Sciences. 2017, No.4 (44), pp. 3–17 (in Russian).

\bibitem{14} V. S. Vladimirov, Methods of the Theory of Generalized Functions, Anal. Methods Special Funct. 6, Taylor and Francis,
London, 2002.

\bibitem{15} Krasnosel’skii M.A., Vainikko G.M., Zabreiko P.P., Rutitcki Ja.B., Stecenko V.Ja., Approximated Solutions of Operator Equations. Walters — Noordhoff, Groningen, 1972.

\bibitem{16} A. N. Tikhonov, A. V. Goncharsky, V. V. Stepanov and A. G. Yagola, Numerical Methods for the Solution of Ill-Posed
Problems, Math. Appl. 328, Kluwer Academic Publishers, Dordrecht, 1995.

\bibitem{17} A. S. Leonov, Solution of Ill-Posed Inverse Problems. Theory Review, Practical Algorithms and MATLAB Demonstrations
(in Russian), Librokom, Moscow, 2010.

\bibitem{18} Engl H.W., Hanke M, Neubauer A. Regularization of Inverse Problems. Dordrecht: Kluwer Academic Publishers,, 1996.

\end{thebibliography}
\end{document}